\documentstyle{amltd}
\begin{document}
\annalsline{159}{2004}
\received{May 25, 2001}
\revised{June 17, 2002}
\startingpage{277}
\def\bye{\end{document}}
 \font\tenrm=cmr10
\def\ritem#1{\item[{\rm #1}]}
\input amssym.tex
\input amssym.def
\input boxedeps.tex 
\SetepsfEPSFSpecial 
\HideDisplacementBoxes
\def\figin#1#2{
$$
 {\BoxedEPSF{#1.eps scaled
#2}%
}%
$$
\noindent}
\catcode`\@=11
\font\twelvemsb=msbm10 scaled 1100
\font\tenmsb=msbm10
\font\ninemsb=msbm10 scaled 800
\newfam\msbfam
\textfont\msbfam=\twelvemsb  \scriptfont\msbfam=\ninemsb
  \scriptscriptfont\msbfam=\ninemsb
\def\msb@{\hexnumber@\msbfam}
\def\Bbb{\relax\ifmmode\let\next\Bbb@\else
 \def\next{\errmessage{Use \string\Bbb\space only in math
mode}}\fi\next}
\def\Bbb@#1{{\Bbb@@{#1}}}
\def\Bbb@@#1{\fam\msbfam#1}
\catcode`\@=12

 \catcode`\@=11
\font\twelveeuf=eufm10 scaled 1100
\font\teneuf=eufm10
\font\nineeuf=eufm7 scaled 1100
\newfam\euffam
\textfont\euffam=\twelveeuf  \scriptfont\euffam=\teneuf
  \scriptscriptfont\euffam=\nineeuf
\def\euf@{\hexnumber@\euffam}
\def\frak{\relax\ifmmode\let\next\frak@\else
 \def\next{\errmessage{Use \string\frak\space only in math
mode}}\fi\next}
\def\frak@#1{{\frak@@{#1}}}
\def\frak@@#1{\fam\euffam#1}
\catcode`\@=12

\newcommand{\R}{{\Bbb R}}
\newcommand{\RN}{{\Bbb R}^n}
\newcommand{\RNC}{\overline{\Bbb R}^n}
\newcommand{\mod}{{\rm mod}\,}
\newcommand{\CAP}{{\rm cap}\,}
\newcommand{\Mod}{{\rm Mod}\,}
\newcommand{\meas}{{\rm meas}\,}
\newcommand{\ar}{{\rm Area}\,}
\newcommand{\g}{\gamma}
\newcommand{\G}{\Gamma}
\newcommand{\va}{\varphi}
\newcommand{\al}{\alpha}
\newcommand{\Arg}{{\rm Arg}\,}
\newcommand{\RC}{{\overline{\Bbb R}}}
\newcommand{\F}{{\cal F}}
\newcommand{\cD}{\cal D}
\newcommand{\e}{\varepsilon}

\newcommand{\be}{\begin{equation}}
\newcommand{\ee}{\end{equation}}
 
\newcommand{\C}{{\Bbb{C}}}
\newcommand{\UH}{{\Bbb{H}}}
\newcommand{\U}{{\Bbb{U}}}
\newcommand{\CC}{\overline{\C}}
\newcommand{\CUH}{\overline{\UH}}
\newcommand{\CU}{\overline{\U}}
\newcommand{\T}{{\Bbb{T}}}
\newcommand{\argu}{{\mbox{arg}\,}}
\newcommand{\om}{\omega}
\newcommand{\D}{\Delta}
\newcommand{\cA}{{\cal{A}}}
\newcommand{\cR}{{\cal{R}}}
\renewcommand{\theequation}{\thesection.\arabic{equation}}

 \title{An isoperimetric inequality\\
for logarithmic capacity of polygons} 
\shorttitle{Logarithmic capacity of polygons} 

 \acknowledgements{This paper was finalized during the first author's visit at the
Technion - Israel Institute of Technology, Spring 2001 under the
financial support of the Lady Devis Fellowship. This author thanks
the Department of Mathematics  of the Technion for wonderful
atmosphere and working conditions during his stay in Haifa. The
research of the first author was supported in part by the Russian
Foundation for Basic Research, grant no.\ 00-01-00118a.}

 \twoauthors{Alexander Yu.\ Solynin}{Victor A.\ Zalgaller}
\institutions{Steklov Institute of Mathematics at St. Petersburg, Russian Academy of
Sciences,\\ St.\ Petersburg, Russia\\
{\eightpoint {\it Current address\/}}: Department of Mathematical Sciences, University of Arkansas,\\ Fayetteville, AR\\
{\eightpoint {\it E-mail address\/}:  solynin@uark.edu}
\vglue6pt
The Weizmann Institute of Science, Rehovot, Israel\\
{\eightpoint {\it E-mail address\/}: victorzalgaller@weizmann.ac.il}}

\centerline{\bf Abstract}
\vglue12pt
We verify an old conjecture of G.~P\'{o}lya and G.~Szeg\H{o}
saying that the regular $n$-gon minimizes the logarithmic capacity
among all $n$-gons with a fixed area. 

\section{Introduction}

\setcounter{equation}{0}

The logarithmic capacity $\CAP E$ of a compact set $E$ in $\R^2$,
which we identify with the complex plane $\C$, is defined by
\be
-\log\CAP E=\lim_{z\to \infty} (g(z,\infty)-\log |z|), \ee where
$g(z,\infty)$ denotes the Green function of a connected component
$\Omega(E)\ni \infty$ of $\CC\setminus E$ having singularity at
$z=\infty$; see \cite[Ch.~7]{G}, \cite[\S11.1]{Po}. By an $n$-gon
with $n\ge 3$ sides we mean a simply connected Jordan domain
$D_n\subset \C$ whose boundary $\partial D_n$ consists of $n$
rectilinear segments called sides of $D_n$. A closed $n$-gon
will be  denoted by $\overline {D}_n$.

Our principal result is 

\specialnumber{1}\proclaim{Theorem} 
For any polygon $D_n$ having a given number of sides $n\ge 3${\rm ,}
\be
\frac{\CAP^2\,\overline{D}_n}{{\rm{Area}}\,D_n}\ge
\frac{\CAP^2\,\overline{D}^*_n}{{\rm{Area}}\,D_n^*}=\frac{n \tan(\pi/
n)\Gamma^2(1+1/n)}{\pi2^{4/n}\Gamma^2(1/2+1/n)} \ee with the sign
of equality only for the regular $n$\/{\rm -}\/gons. \endproclaim 

In Theorem~1 and
below, $\Gamma(\cdot)$ denotes the Euler gamma function and
 $D^*_n$ stands for the regular $n$-gon centered at $z=0$
with one vertex at $z=1$.

In other words, Theorem~1 asserts that the regular closed polygon has the
minimal logarithmic capacity among all closed polygons with a fixed
number of sides and prescribed area. For $n\ge 5$, this solves an
old problem  posed by G.~P\'{o}lya and G.~Szeg\H{o} \cite{PS}. For
$n=3,4$, the problem was solved by P\'{o}lya and Szeg\H{o}
themselves \cite[p.158]{PS}. Their method based on Steiner
symmetrization allows them to establish similar isoperimetric
inequalities for the conformal radius, torsional rigidity,
principal frequency, etc. However it fails for $n\ge 5$ since
Steiner symmetrization increases dimension (= number of sides) of
a polygon in general. In \cite[p.159]{PS} the authors note that
``to prove  (or disprove) the analogous theorems for  regular
polygons with more than four sides is a challenging task''.

For the conformal radius this task was solved in \cite{So1},
where it was shown that the regular $n$-gon maximizes the conformal
radius among all polygons with a given number $n\ge 3$ of sides and with
a prescribed area. The
present  work proves the P\'{o}lya-Szeg\H{o} conjecture for the
logarithmic capacity. For the torsional rigidity and principal
frequency the problem is still open.

A similar question concerning the minimal logarithmic capacity
among all compact sets with a prescribed perimeter is nontrivial
only for convex sets. This question was studied by G.~P\'{o}lya
and M.~Schiffer and Chr.~Pommerenke, see \cite[p.~51, Prob.~11]{Po},
who proved that a needle (rectilinear segment) is a unique
minimal configuration of the problem. Since a needle can be viewed
as a degenerate $n$-gon, there is no difference between the
convex polygonal case
and the general case. Thus the regular $n$-gons do not
minimize the logarithmic
capacity over the set of all $n$-gons with a prescribed perimeter.
To the contrary, they provide the maximal value for this problem; see
\cite[Th.~10]{So1.1}.

Any isoperimetric problem for polygons of a fixed dimension can be
considered as a discrete version of an isoperimetric problem among
all simply connected (or more general) domains. It is interesting
to note that solutions to continuous versions for the above
mentioned functionals have been known for a long time; cf.~\cite{PS}.
The discrete problems are much harder. The situation here is
opposite to the classical isoperimetric area-perimeter problem,
where solution to the continuous version requires much stronger
techniques than the discrete case.

The idea of the proof in \cite{So1}, used also in the present
paper, traces back to the classical  method of finding the
area of a polygon:  divide a polygon into triangles and use the
additivity property of the area. Although the characteristics
under consideration are not additive functions of a set, often
they admit a certain kind of ``semiadditivity'', at least for
special decompositions. For instance, the reduced module
$m(D,z_0)$ of a polygon $D$ at its point $z_0\in D$,  a
characteristic linked with the conformal radius and logarithmic
capacity,  admits an explicit upper bound $B$ given by a weighted
sum of the reduced modules of triangles composing $D$, each of
which has a distinguished vertex at~$z_0$. The precise definitions
and formulations will be given in Section~2. This explicit bound
$B$ is a complicated combination of functions including the Euler
gamma function, which depends on the angles and areas of triangles
composing $D$. For the problem on the conformal radius, it was
shown in \cite{So1} that the corresponding maximum of $B$ taken
among all admissible values of the parameters provides the sharp
upper bound for the reduced module $m(D,z_0)$ where $\ar
D$ is fixed.

For the logarithmic capacity when the same method is applied, the
situation is different; the explicit upper bound $B$ contains more
parameters and the supremum of $B$ among all admissible
decompositions of $D$ into triangles is infinite. Even more, for
instance for the regular $n$-gon there is only one
decomposition (into equal triangles) that gives the desired upper
bound for the reduced module. All other decompositions lead to a
bigger upper bound and therefore should be excluded from
consideration if we are looking for a sharp result.

So it is important to select a more narrow subclass of
decompositions among which the maximal value of $B$ corresponding
to the logarithmic\break capacity is finite and provides the sharp bound
for the considered characteristic of $D$. This is the subject of
our study in Section~3. The selected subclass contains
decompositions of $D$ into triangles that are proportional in a
certain sense. This result is of independent interest. We present
it in our Theorem~2 restricting for simplicity of formulation to
the case of convex polygons. The general version for the
nonconvex case is given by Theorem~4 in Section~3.

Let $D_n$ be a convex $n$-gon having vertices
$A_1,\ldots,A_n,A_{n+1}=A_1$ enumerated in the positive direction
on $\partial D_n$. A system of Euclidean triangles
$\{T_k\}_{k=1}^n$ is called admissible for $D_n$ if $T_k\cap
D_n\not=\emptyset$, $T_k$ has the segment $[A_k,A_{k+1}]$ as its
base, and if for all $k=1,\ldots,n$, $T_k$ and $T_{k+1}$ have a
common boundary segment which is an entire side of at least one of
these triangles but not necessarily of both of them.

In Section~3, we give a more general definition of admissibility
for a system of triangles suitable for nonconvex polygons. For
a convex polygon, the definition of admissibility presented above
and the definition given in Section~3 are equivalent.

Let $\al_k$ denote the angle of $T_k$ opposite the base
$[A_k,A_{k+1}]$. An admissible system $\{T_k\}_{k=1}^n$ is called
{\it proportional} if the quotient $\al_k/\ar T_k$ does not depend
on $k=1,\ldots,n$.

\specialnumber{2}\proclaim{Theorem} 
For every convex $n$-gon $D_n$ there is at least one
proportional system $\{T_k\}_{k=1}^n$ that covers $D_n${\rm ,} i.e.
\be
\bigcup_{k=1}^n \overline{T}_k\supset {D}_n.
\ee
\endproclaim

 Theorem~2 is sharp in the sense that there are polygons, for
instance, triangles and regular $n$-gons, that have a unique
proportional system satisfying (1.3). For triangles, Theorem~2
provides a good exercise for the course of elementary geometry. It
is not difficult to show that any rectangle different from a
square admits a parametric family of proportional systems
satisfying (1.3).
Figures~1a)--1c) show possible types of proportional configurations
for a rectangle $\cR$: a)~a~proportional system that does not cover $\cR$;
b)~a~proportional covering system consisting of disjoint triangles;
c)~a~proportional covering system consisting of overlapping triangles
the union of which
is strictly larger than $\cR$ (if $\cR$ is sufficiently long).
Figure~1d), which is a slightly modified version of Figure~1c), gives
an example of a proportional system of six triangles for a nonconvex
hexagon. As we have already mentioned,
the precise definitions for the nonconvex case will be given in Section~3.
\vglue6pt
\centerline{\BoxedEPSF{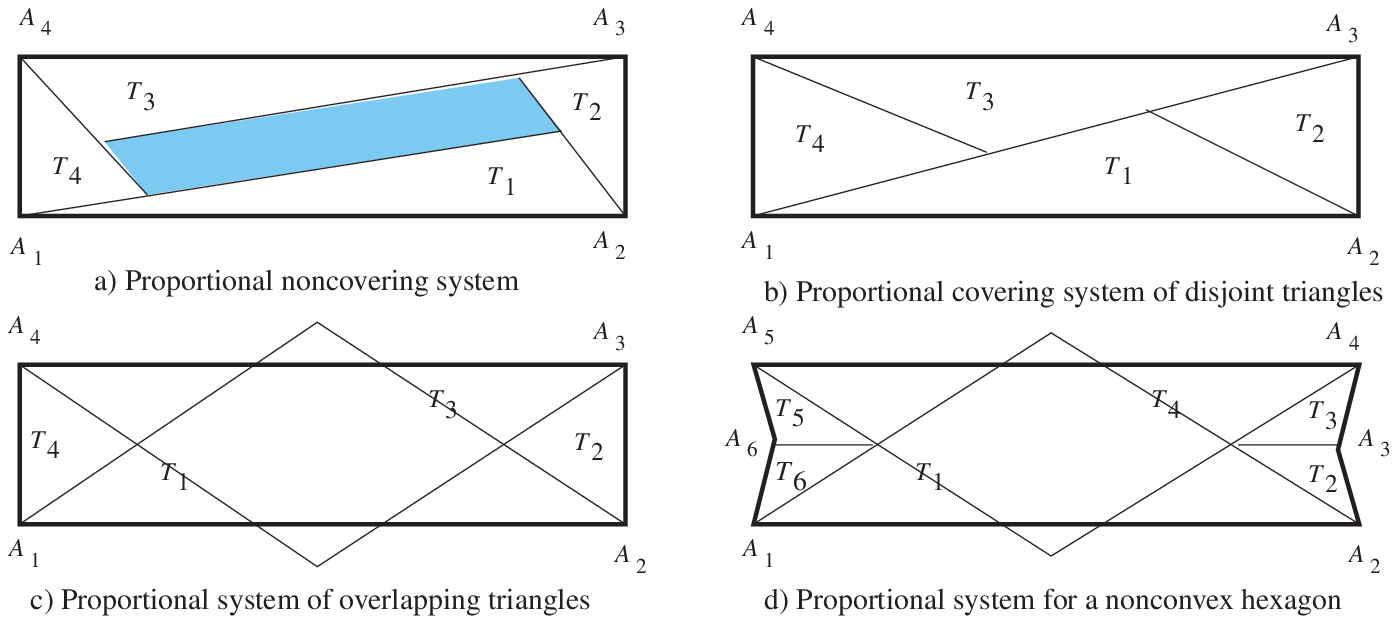 scaled 900}}

 \medbreak \centerline{Figure 1. Proportional systems of triangles}
\vglue6pt

To prove a generalization of Theorem~2 for the nonconvex case, we
show in Lemma~5 that the family of all proportional systems for
$D_n$ admits a natural continuous parametrization. Then the  continuity
property is used in Lemma~6 to show that at least one system of
any continuously parametrized family of admissible systems covers
$D_n$. It is important to note that Theorems~2 and 4 possess
counterparts in other cases of proportionality between some two
characteristics of a triangle (not necessarily the base angle and
the area).

Section~4 finishes the proof of Theorem~1.

The subject of this paper lies at the junction of potential
theory, analysis, and geometry. And this work is a natural result
of combined efforts of an analyst and a geometer.

We are grateful to the referees for their constructive criticism
and many valuable suggestions, which allow us to improve the
exposition of our results. In particular, the short proof of
Lemma~2 in Section~4 was suggested by one of the referees.
 
\section{Logarithmic capacity and reduced module}
 
There are several other approaches to the measure of a set
described by the logarithmic capacity. For example, the geometric
concept of transfinite diameter due to M.~Fekete and the concept
of the Chebyshev's constant from polynomial approximation lead to
the same characteristic; cf.\ \cite[\S~10.2]{Du},\break \cite[Ch.~7]{G}.

If a compact set $E$ is connected, then $\Omega(E)$ is a simply
connected domain containing the point at $\infty$. In this case
the logarithmic capacity is equal to the {\it outer radius} $R(E)$
defined as follows. Let  $$f(z)=z+a_0+a_1 z^{-1}+\ldots$$ map
$\Omega(E)$ conformally onto $|\zeta|>R$. The radius $R=R(E)$ of
the omitted disk is uniquely determined and is called the outer
radius of $E$; see \cite[\S~10.2]{Du}, \cite[Ch.~7]{G}.

The outer radius $R(E)$ can be considered as a characteristic of a
simply connected domain $\Omega(E)$ at its point at $\infty$.
Another approach due to\break O.~Teichm\"{u}ller leads to essentially
the same characteristic of a simply connected domain.  For $R>0$
big enough, let $\Omega_R(E)$ be a doubly connected domain between
$E$ and the circle $C_R=\{z:\,|z|=R\}$ and let $\mod
(\Omega_R(E))$ denote the module of $\Omega_R(E)$ with respect to  the
family of curves separating the boundary components of
$\Omega_R(E)$; see \cite[Ch.~2]{J}. Then there is a finite limit
\be
m(\Omega(E),\infty)=\lim_{R\to
\infty}(\mod(\Omega_R(E))-(1/2\pi)\log R) \ee called {\it the
reduced module of $\Omega(E)$ at $z=\infty$}. The reduced module
can be defined for any point $a\in \Omega(E)$ finite or not;
cf.\cite[Ch.~2]{J}  but we shall use this notion with $a=\infty$
only. It is well known \cite[\S~1.3]{D}, \cite[Ch.~7]{G} that
\be
m(\Omega(E),\infty)=-(1/2\pi)\log \CAP E. \ee Thus, (1.2) holds if
and only if $\Omega(\overline{D}^*_n)$ has the maximal reduced module at
$\infty$ among all domains $\Omega(\overline{D}_n)$ corresponding to polygons
$D_n$ such that $\ar D_n=\ar D^*_n$.

As   mentioned in the introduction, to prove Theorem~1 we apply
the method developed in \cite{So1}, \cite{So2} based on a special
triangulation of $\Omega(E)$.

By a trilateral $D=D(a_0,a_1,a_2)$ we mean a simply connected
domain $D\subset \CC$ having three distinguished points $a_0$,
$a_1$, and $a_2$ called {\it vertices} on its boundary. Each
trilateral will have a distinguished side called {\it the base};
the opposite vertex and angle will be called the base vertex and
the base angle respectively. For our purposes it is enough to deal
with trilaterals having the vertex $a_0$ at $\infty$ with a
piecewise smooth Jordan boundary such that $l_R=D\cap C_R$
contains   only one connected component for all $R>0$ sufficiently
large. Let $D_R=D\cap \U_R$, where $\U_R=\{z:\,|z|<R\}$.
Considering $D_R$ as a quadrilateral with distinguished sides
$\widehat{a_1a_2}$ and $l_R$, let $\mod (D_R)$ denote the module
of $D_R$ with respect to the family of curves separating $\widehat{a_1a_2}$
from $l_R$ in $D_R$; cf.\break \cite[Ch.~2]{J}. Let $D$ have an inner
angle $0<\varphi\le 2\pi$ at $a_0=\infty$. The limit
\be
m(D;\infty|a_1,a_2)=\lim _{R\to \infty}(\mod (D_R)
-(1/\varphi)\log R), \ee provided that it exists and is finite, is
called {\it the reduced module of $D$ at $a_0=\infty$}. This
notion was introduced in \cite{So1}. In \cite{So3} some sufficient
conditions for the existence of the limit in (2.3) are given. In
this paper we deal with rectilinear trilaterals only which
guarantees existence of all the reduced modules considered below.

Regarding the infinite circular sector $P=P(\rho,\alpha)=
\{z:\,|z|>\rho, 0<\arg
z<\al\}$, $\rho>0$, $0<\al\le 2\pi$ and the upper half-plane $\UH=\{z:\,\Im
z>0\}$ as trilaterals with vertices $\infty$, $\rho$, $\rho e^{i\al}$ and
$\infty$, $0$, $\rho$, respectively, and computing the corresponding
limits in (2.3), we get,
\be  
 m(P;\infty|\rho, \rho e^{i\al})=-(1/\al)\log \rho, \quad 
m(\UH;\infty|0,\rho)= (1/\pi)\log (4/\rho), \quad\enspace
\ee
which provides two
useful examples of the reduced modules.

 The change in the reduced module under conformal
mapping can be worked out by means of a standard formula
\cite{So1}, \cite{So3}: if a function $f(\zeta)=A\zeta^\al(1+o(1))$ with
$\al>0$, $A\not=0$, and $o(1)\to 0$ as $\zeta\to \infty$ maps the upper
half-plane $\UH$ conformally onto a trilateral $D=D(a_0,a_1,a_2)$,
$a_0=\infty$ such that $f(\infty)=\infty$, $f(0)=a_1$, $f(1)=a_2$,
then
\be
m(D;\infty|a_1,a_2)=(1/\pi)\log 4 -(1/(\alpha \pi))\log |A|.
\ee

Let $T_1,\ldots,T_n$ be pairwise disjoint trilaterals in a simply
connected domain~$D$, $\infty\in D\subset \CC$, such that $T_k$
has a vertex $a_0^k$ at $\infty$ and the opposite side
$\widehat{a_1^ka_2^k}$ on $\partial D$;  see Figure~2, where for simplicity
the point at $\infty$ is represented by a finite point $a_0$. The
next result from \cite{So1} linking the reduced module of $D$ with
the reduced modules of trilaterals of its decomposition, is basic
for our further considerations.

\proclaimtitle{\cite{So1}}
\specialnumber{3}\proclaim{Theorem}  Let $T_k$ have an angle $0<2\pi \al_k<2\pi$ at the
vertex $a_0^k$ and  for every $k=1,\ldots,n$ let the reduced module
of $T_k$ at $\infty$ exist. If $\sum_{k=1}^n \al_k=1${\rm ,} then
\be
m(D,\infty)\le \sum_{k=1}^n \al^2_k m(T_k;\infty|a_1^k,a_2^k). \ee
Let $f$ map $D$ conformally onto $\U^*=\CC\setminus \CU_1$ such
that $f(\infty)=\infty$. Equality occurs in {\rm (2.6)} if and only if
for every $k=1,\ldots,n${\rm ,}  $f(T_k)$ is an infinite circular sector
{\rm (}\/of opening $2\pi\al_k${\rm )}  and if the vertices of $T_k$ correspond
under the mapping $f$ to the geometric vertices of this sector.
\endproclaim

\figin{fig4}{700}
\centerline{Figure 2. Decomposition into
trilaterals} 
\vglue12pt

 The proof of (2.6) in \cite{So1} is based on basic properties of
the extremal length. Another approach to more general problems on
the extremal decomposition developed by V.~N.~Dubinin \cite{D}
uses the theory of capacities.

 Now we consider an instructive example that is important
for what then follows. Up to the end of the paper all considered
trilaterals will be rectilinear triangles (finite or not)
having their geometric
vertices as the distinguished boundary points. In this case we shall use
the terms ``triangle'' and ``infinite triangle'' instead of ``trilateral''.
Thus, everywhere below, ``triangle'' means a usual Euclidean triangle.

For $\al>0$, $\beta_1>0$, $\beta_2>0$ such that
$$
\beta_1+\beta_2=1+2\al,
$$
and $a>0$, let $T=T(\al,\beta_1,a)$ be the triangle having
vertices at $a_0=0$, $a_1=a$, and $a_2=e^{i2\pi\al}(a \sin
\pi\beta_1/\sin \pi\beta_2) $ and the side $[a_1,a_2]$ as its
base. Then $T$ has  interior angles $2\pi\al$, $\pi(1-\beta_1)$,
and $\pi(1-\beta_2)$ at the vertices $a_0$, $a_1$, and $a_2$,
respectively. Let $V_\al=\{z:\,0<\arg z <2\pi\al\}$, and let
$S(\al,\beta_1,a)=V_\al \setminus \overline{T}$. Then
$S=S(\al,\beta_1,a)$ is an infinite rectilinear triangle having
vertices at\break $a_\infty=\infty$, $a_1$, and $a_2$, which will be
called {\it the sector associated with} $T$. In Section~3, the
notion of the associated sector will be used in a more general
context.

To find the reduced module $m(S;\infty|a_1,a_2)$, we consider
the Schwarz-Christoffel function
\be
f(\zeta)=a-e^{-i\pi\beta_1}\,C\,\left(\int_0^\zeta t^{\beta_2-1}
(1-t)^{\beta_1-1}\,dt
-B(\beta_1,\beta_2)\right)
\ee
with
$$
C=\frac{a \sin 2\pi\al}{\sin
\pi\beta_1 B(\beta_1,\beta_2)},
$$
where $B(\cdot,\cdot)$
denotes the Euler beta function. The function $f$ maps the upper
half-plane $\UH$ conformally onto the infinite triangle $S$
such that $f(\infty)=\infty$, $f(1)=a_1$, $f(0)=a_2$.
From (2.7),
\be
f(\zeta)=(C/2\al)\zeta^{2\al} +\hbox{ \rm constant } +o(1),
\ee
where $o(1)\to 0$ as $\zeta\to \infty$.

From (2.5), (2.7), and (2.8), using the second equality in (2.4) with $\rho=1$,
we obtain the desired
formula for the reduced module of $S$:
\be
m(S;\infty|a_1,a_2)=\frac{1}{2\pi\al}\log\frac{2^{4\al+1} \al
B(\beta_1,\beta_2)\, \sin \pi\beta_2}{a\sin2\pi\al}.
\ee

Let 
$s=\ar T$ be the area of the triangle $T=T(\al,\beta_1,a)$.  Then
from elementary trigonometry, $$a=\left[\frac{2s\sin
\pi\beta_2}{\sin2\pi\al \sin\pi\beta_1}\right]^{1/2}.$$
Substituting this in (2.9), we get
\be
m(S;\infty|a_1,a_2)=\frac{1}{2\pi\al}\log\frac{2^{4\al+1}\al
B(\beta_1,\beta_2)
(\sin\pi\beta_1 \sin\pi\beta_2)^{1/2}}{(2s\sin2\pi\al)^{1/2}}.\qquad
\ee

For a fixed $\al$, $0<\al <1/2$, and $s>0$, let $F(\beta_1)$
denote the right-hand side of (2.10) with $\beta_2=1+2\al -\beta_1$
regarded as
a function of $\beta_1$, $2\al<\beta_1 <1$. The next lemma shows that
$F$ is concave in
$2\al<\beta_1<1$. This implies, in particular, that the isosceles
infinite triangle $S(\al,1/2+\al,a)$ has the maximal reduced module among
all infinite triangles $S(\al,\beta_1,a)$ with  fixed angle $2\pi\al$ and
fixed area $s$ of $T(\al,\beta_1,a)$.

 \specialnumber{1}\proclaim{Lemma}  Let $0<\al<1/2$ and $s>0$ be fixed. Then $F(\beta_1)$ is
strictly concave in $2\al<\beta<1$ and satisfies the equation
$F(\beta_1)=F(\beta_2)$ for $2\al<\beta_1\break <1$. In particular{\rm ,} \be
F(\beta_1)<F(1/2+\al)=\frac{1}{2\pi\al}\log\frac{4^\al \al
B(1/2,1/2+\al)}{(s\tan \pi\al)^{1/2}} \ee for $2\al<\beta_1<1$
such that $\beta_1\not= 1/2+\al$. \endproclaim

\demo{Proof} Since $B(\beta_1,\beta_2)=\Gamma(\beta_1)\Gamma(\beta_2)/
\Gamma(\beta_1+\beta_2)$ and $\beta_1+\beta_2=1+2\al$,
(2.10) implies
\be
F(\beta_1)=
\frac{1}{2\pi\al}\log\frac{2^{4\al +1} \al
\G(\beta_1)\Gamma(\beta_2) \left(\sin \pi\beta_1\,\sin
\pi\beta_2\right)^{1/2}}{\G(1+2\al) \left(2s\sin
2\pi\al\right)^{1/2}}. \qquad\quad\ee Using the reflection formula
$$
\G(z)\G(1-z)=\pi/\sin \pi z,
$$
from (2.12) we obtain
$$
F(\beta_1)=\frac{1}{4\pi\al}\log \frac{\G(\beta_1)\G(\beta_2)}{\G(\beta_1-2\al)
\G(\beta_2-2\al)}+\frac{1}{2\pi\al}\log
\frac{2^{4\al+1}\pi\al}{\G(1+2\al)\left(2s\sin 2\pi\al\right)^{1/2}},
$$
where $\beta_2=1+2\al-\beta_1$ and the second term does not
depend on $\beta_1$.
Differentiating twice, we find
\be
F''(\beta_1)=\frac{1}{4\pi\al}\left[ \psi'(\beta_1)-\psi'(\beta_1-2\al)+
\psi'(\beta_2)-\psi'(\beta_2-2\al)\right]<0,\quad
\ee
which is negative because $\psi'(z)=\sum_{k=0}^\infty (t+k)^{-2}$ strictly
decreases for $t>0$ (\cite[p.~45]{BE}). Here and below,
$\psi$ denotes the logarithmic derivative
of the Euler gamma function.
Inequality (2.13) shows that $F(\beta_1)$ is strictly concave.

Since $\beta_2=1+2\al -\beta_1$, the symmetry formula
$F(\beta_1)=F(1+2\al-\beta_1)$ follows immediately from (2.10).
Symmetry and concavity properties  imply that $F$ takes its
maximal value at $\beta_1=1/2+\al$. Substituting $\beta_1=1/2+\al$
in (2.10) and using the formula
$B(1/2+\al,1/2+\al)=4^{-\al}B(1/2,1/2+\al)$, we get (2.11), and the
lemma follows. \enddemo

Let $S_n=S(1/n,1/2+1/n,a)$  with  $a=(2s/\sin(2\pi/n))^{1/2}$.
Then (2.11) with $\al=1/n$, $\beta_1=1/2+1/n$ gives
$$
m(S_n;\infty|a_1,a_2)=\frac{n}{4\pi} \log\frac{\pi
4^{2/n}\Gamma^2(1/2+1/n)}{s\tan(\pi/n)\Gamma^2(1/n)}.
$$

The latter relation combined with the assertion on the equality cases
in Theorem~3 leads to the well-known formula for the
reduced module of the exterior of the regular $n$-gon $D^*_n(A)$
having the area $A$; cf.\ \cite[p.273]{PS}:
\be
m(\Omega(\overline{D}^*_n(A)),\infty)=\frac{1}{4\pi}\log\frac{\pi}{A}
\frac{4^{2/n}n\Gamma^2(1/2+1/n)}{\Gamma^2(1/n) \tan(\pi/n)}. \ee

The next lemma treats $m(\Omega(\overline{D}^*_n(A)),\infty)$ as a function
of the number of sides of $D^*_n(A)$.

 \specialnumber{2}\proclaim{Lemma}  For a fixed area $A${\rm ,} the reduced
module $m(\Omega(\overline{D}^*_n(A)),\infty)$ is strictly increasing in
$n$.\endproclaim

Lemma~2 easily follows from the concavity result of Lemma~7, and
the proof  is given in \pagebreak Section~4.

\section{Triangular covers of a polygon}
 
 To prove Theorem~1 for the nonconvex
polygons, we need a generalization of Theorem~2 for this case.
First we fix terminology and necessary notation.
 Let $T$ be a triangle with the base
$[a_1,a_2]$ and the base vertex $a_0$. Let $V$ be the smallest
infinite sector with the vertex at $a_0$ that contains $T$. Then the
infinite triangle $S=V\setminus \overline {T}$ having the base
$[a_1,a_2]$ and the base vertex $a_\infty=\infty$ will be called
{\it the sector associated with} $T$.

Let $D$ be an $n$-gon with vertices $A_1,\ldots,A_n$. Throughout
this section we use the following conventions concerning
numbering:
\begin{itemize}
\item[i)] Cyclic convention: if a system $\{x_k\}_{k=1}^n$ contains $n\ge
1$ elements, then $x_{n+1}:=x_1$, $x_0:=x_n$, etc.

\item[ii)] Positive orientation convention: numeration of geometric
objects, e.g. vertices, angles, sides of a polygon $D$, triangles
covering $D$, etc. agrees with the positive orientation on
$\partial D$.
\end{itemize}

A triangle $T$ having an associated sector $S$ is called {\it
admissible} for $D$ if $T\cap D\not=\emptyset$, the base of $T$
lies on $\partial D$ (the base of $T$ need not consist of an
entire side of $D$), $S\cap D=\emptyset$, and if each (closed)
side (not base!) of $S$ contains at least one vertex of $D$. Of
course, the first condition follows from the second and third
conditions.

A system of triangles $\{T_i\}_{i=1}^m$ is called {\it admissible}
for $D$ if each $T_i$ is admissible, the associated sectors $S_i$
are pairwise disjoint and if $\mathbold{\cup}_{i=1}^m\overline {S}_i$ covers
the complement of
the convex hull $\hat D$ of $D$.

If $T_i$ has the base angle $\al_i$, which is equal to the angle of $S_i$
at $z=\infty$, the latter conditions imply that
$\sum_{i=1}^m \al_i=2\pi$.

It is important to emphasize that for the case of convex polygons this
definition of  admissibility is equivalent to the definition
of the  admissibility
of a system of triangles given in the introduction.

An admissible system $\{T_i\}_{i=1}^m$ is called {\it regular} if
for $i=1,\ldots,m$, each side of $S_i$ contains only one
vertex of $D$.

Let $\al_i$ and $\sigma_i$ denote the base angle and area of
$T_i$. By the {\it coefficient} of $T_i$ we mean the quotient
$k_i=\al_i/\sigma_i$. An admissible system $\{T_i\}_{i=1}^m$ is
called {\it proportional} if $k_1=\ldots =k_m$.

In this terminology, the system of triangles shown in Figure~1d)
is admissible,
regular, and proportional, which covers the hexagon
for which it is constructed.

The purpose of this section is to prove the following
theorem, which includes Theorem~2 as a special case.
 
\specialnumber{4}\proclaim{Theorem}  For every $n$\/{\rm -}\/gon $D$ whose convex hull $\hat D$ has $\hat n\ge 3$
sides{\rm ,} there is at least one proportional
system $\{T_i\}_{i=1}^m$, $\hat n \le m\le n${\rm ,}  that
covers $D${\rm ,}
i.e.
\be
\bigcup_{i=1}^m \overline T_i\supset  D. \ee In particular{\rm ,}
 $$
\sum_{i=1}^m \ar T_i\ge \ar D. $$ \endproclaim 

The proof of Theorem~4 will be given after Lemmas~5 and 6 which
study the family of all proportional systems $\{T_i\}$ admissible
for  a given $n$-gon $D$ in general position. The latter means
that no three vertices of $D$ belong to the same straight line and
no side or diagonal is parallel to any other side or diagonal.
For such $D$, we show in Lemma~5  that the set of all proportional
systems admits a natural continuous parametrization.

To prove Lemma~5, we need the following variant of the standard
implicit function theorem.

\specialnumber{3}\proclaim{Lemma}  Let $u_i${\rm ,} $i=1,\ldots,n+1${\rm ,} be real\/{\rm -}\/valued functions having
continuous partial derivatives in a neighborhood of $x^0\in
\R^{n+1}$. Let $u_i=u_i(x_{i-1},x_i,x_{i+1})${\rm ,} $i=2,\ldots,n+1${\rm ,}
$u_1=u_1(x_1,x_2)$ and suppose that the $u_i$ do not depend on the other variables. If
for $x=x^0${\rm ,}
\be
u_1(x)=u_2(x)=\ldots=u_{n+1}(x) \ee and if the partial derivatives
$u_{i,j}$ satisfy \be u_{i,i+1}(x^0)>0,\quad u_{i,j}(x^0)\le 0 \ee
for all $i$ and $j\not=i+1${\rm ,} then for every $x_1$ in some small
interval $(x_1^0-\delta,x_1^0+\delta)$ equations {\rm (3.2)} define a
unique solution $x_2=x_2(x_1),\ldots,x_{n+1}=x_{n+1}(x_1)$
continuously differentiable in $(x_1^0-\delta,x_1^0+\delta)$ and
such that $x_i(x_1^0)=x_i^0$.

Let the functions $u_i$ satisfy the following additional
assumptions\/{\rm :}\/
\begin{itemize}
\ritem{1)} For every $u_i(x_{i-1},x_i,x_{i+1})$ depending on three
parameters{\rm ,} the neighbors $u_{i-1}(x_{i-1},x_i)$ and
$u_{i+1}(x_{i+1},x_{i+2})$ depend on two parameters each and
$u_{i-1,i-1}(x^0)=0${\rm ,} $u_{i,i-1}(x^0)<0$.

\ritem{2)} If $u_i(x_i,x_{i+1})$ and $u_{i+1}(x_{i+1},x_{i+2})$ depend on
two parameters each{\rm ,} then $u_{i+1,i+1}(x^0)<0$.
\end{itemize}

Then
\be
x'_j(x_1)>0 \ee for every index $j$ such that
$u_j=u_j(x_j,x_{j+1})$ depends on two parameters and every $x_1\in
(x_1^0-\delta,x_1^0+\delta)$. \endproclaim

\demo{Proof} Setting $v_i:=u_{i+1}-u_1$, we consider the equations
\be
v_i(x_1,\ldots,x_{n+1})=0, \quad i=1,\ldots,n. \ee The assumptions
of the lemma imply that the Jacobian $\Delta_n=|v_{i,j}(x^0)|$ in
variables $x_2,\ldots,x_{n+1}$ has the following form
\be
\Delta_n=\left|
\begin{array}{cccccccccc}-&+& & & & & & & & \\
-&\div&+& & & & & & & \\-&\div&\div&+& & & & & & \\\cdot&\cdot
&\cdot&\cdot&\cdot& & & & & \\ \cdot&\cdot &\cdot
&\cdot&\cdot&\cdot& & & &
\\\cdot&\cdot &\cdot &\cdot &\cdot&\cdot&\cdot& & &
\\\cdot&\cdot &\cdot &\cdot &\cdot &\cdot&\cdot&\cdot& &   \\
-&\div &\cdot &\cdot &\cdot &\cdot &\div&\div&+& \\-&\div &\cdot
&\cdot &\cdot &\cdot &\cdot &\div&\div&+\\-&\div &\cdot &\cdot
&\cdot &\cdot &\cdot &\cdot &\div&\div\\\end{array}\right| \ee
Here $+$, $-$, and $\div$ represent  positive, negative, and
nonpositive elements, respectively. The blank spaces are supposed
to be filled with zeros.

We claim that
\be
\Delta_n\not=0 \quad{\mbox{and}}\quad \Delta_n/|\Delta_n|=(-1)^n.
\ee The proof is by induction. For $n=1,2,3$ the result is
obvious. Assume the assumption holds true for the dimension $n-1$.
Then expanding $\Delta_n$ in its last column, we get
\be
\Delta_n=c_{n,n}\Delta'_{n-1}-c_{n-1,n}\Delta''_{n-1}, \ee where
$c_{n,n}=v_{n,n+1}(x^0)\le 0$, $c_{n-1,n}=v_{n-1,n+1}(x^0)>0$ and
$\Delta'_{n-1}$, $\Delta''_{n-1}$ are  $(n-1)$-dimensional
determinants of the form (3.6). The inductive assumption and  (3.8)
imply (3.7), which
 shows that $v_1,\ldots,v_n$ satisfy the assumption of the
standard implicit function theorem. Therefore (3.5), or
equivalently (3.2), defines in $(x_1^0-\delta,x_1^0+\delta)$ the
unique continuously differentiable functions $x_i=x_i(x_1)$,
$i=2,\ldots,n+1$.

Now we shall prove the additional assertion (3.4). To check the
sign of $y_{i-1}=x'_i(x^0_1)$, we note that $y_1,\ldots,y_n$
satisfy the system of linear equations
\be
(v_{i,j}(x^0))\,(y_i)=(-v_{i,1}(x^0)), \ee where $(v_{i,j}(x^0))$
is the Jacobi matrix corresponding $\Delta_n$.  By Cramer's
rule, \be y_k=\Delta^k_n/\Delta_n, \ee where the determinant
$\Delta^k_n$ is obtained from $\Delta_n$ by replacing its $k$-th
column with the column $(-v_{i,1}(x^0))$.

\pagebreak

1) If $u_2=u_2(x_2,x_3)$ depends on two parameters, then \begin{eqnarray}
v_{1,1}(x^0)&=&\ldots=v_{n-1,1}(x^0)=-u_{1,1}(x^0)>0,\\
v_{n,1}(x^0)&=&u_{n+1,1}(x^0)-u_{1,1}(x^0)>0.\nonumber
\end{eqnarray} This shows that the
determinant $\Delta^1_n$ has the form (3.6). Therefore
$\Delta^1_n\not=0$, $\Delta^1_n/|\Delta^1_n|=(-1)^n$. This
combined with (3.7) and (3.10) implies that $y_1=$\break \vglue-10pt \noindent $x'_2(x_1^0)>0$.
\vglue6pt
2) If $u_2$ depends on three parameters, then by the assumptions
of the lemma, $u_3=u_3(x_3,x_4)$ depends on two parameters and
\begin{eqnarray} 
v_{1,1}(x^0)&=&u_{2,1}(x^0)-u_{1,1}(x^0)=u_{2,1}(x^0)<0,\\
v_{2,1}(x^0)&=&\ldots=v_{n-1,1}(x^0)=0,\quad
v_{n,1}(x^0)=u_{n+1,1}(x^0)>0.\nonumber \end{eqnarray} Therefore,
$\Delta^2_n$ has the form (3.6) in the case under consideration.
Hence, $x'_3(x^0_1)=\Delta^2_n/\Delta_n>0$.

Note that the assumptions of the lemma allow us to apply the same
arguments for the ``shifted'' functions: \be v_i^2:=u_{i+1}-u_2,
\quad i=2,\ldots,n+1, \ee if $u_2$ satisfies additional assumption
1), or for the functions \be v_i^3:=u_{i+1}-u_3, \quad
i=3,\ldots,n+2,\ee if $u_2$ and $u_3$ satisfy additional
assumption 2). Since (3.13) and (3.14) are equivalent to (3.5),
each of them defines the same system of solutions
$x_1,x_2,\ldots,x_{n+1}$ in a neighborhood of $x^0$.
Considerations above show that $x_3'(x_1)=x_3'(x_2)x'_2(x_1)\break >0$ in
the case corresponding to (3.13) and
$x'_4(x_1)=x_4'(x_2)x_2'(x_1)>0$ in the case corresponding to
(3.14).

Repeating these arguments, after a finite number of steps we get the desired assertion (3.4). \enddemo

The next geometrically obvious lemma will be used in the proofs of\break
Lemma~5 and Theorem~4. Let $T$ be a triangle with the base
$[a_1,a_2]$. A system of triangles
$\{T_i\}_{i=1}^m$ is called admissible for $T$ if:
\begin{itemize}
\item[1)] $T_i$ has a base $[a_{1,i},a_{2,i}]$, such that the segments
$[a_{1,i},a_{2,i}]$, $i=1,\ldots,m$, constitute a disjoint
decomposition of the base $[a_1,a_2]$;

\item[2)] $T_i$ and $T_{i+1}$ are disjoint and  have a common boundary segment
that is an entire side
of at least one of these triangles but not necessarily of both of them;

\item[3)] $T_1$ and $T_m$ each has a common boundary segment with
$\partial T\setminus [a_1,a_2]$.
\end{itemize}
\specialnumber{4}\proclaim{Lemma} 
If $\{T_i\}_{i=1}^m$ is admissible for $T$, then $\mathbold{\cup}_{i=1}^m
\overline{T}_i \supset T$.
\endproclaim
 
It is important to emphasize that all the triangles under consideration
are the  usual Euclidean triangles. An elementary inductive proof of
the lemma is left to the readers.

Let $D$ be an $n$-gon in  general position with vertices
$A_1,\ldots,A_n$. Let the convex hull $\hat D$ of $D$ have
vertices $A_1'=A_1,A_2',\ldots, A'_{\hat n}$. If $A_1'=0$ and
$A_2'>0$, we say that $D$ is in standard position. Let
$\{T_i\}_{i=1}^m$ be an admissible system for $D$. Let
$[a_{1,i},a_{2,i}]$ and $a_{0,i}$ denote the base and base vertex
of $T_i$. Let $\gamma_{1,i}$ and $\gamma_{2,i}$ be the closed
sides of the associated sector $S_i$ starting at the points
$a_{1,i}$ and $a_{2,i}$, respectively. Since $D$ is in general
position, $\gamma_{k,i}$ contains one or two vertices of $D$. We
shall denote them by $B_{k,i}$ and $B_{k,i}'$, where the second
vertex, if it exists, lies between $B_{k,i}$ and $a_{k,i}$. Let
$l_{k,i}$ and $l_{k,i}'$ denote the rays outgoing from $B_{k,i}$
and $B'_{k,i}$ each containing the side $[a_{k,i},a_{0,i}]$ of
$T_i$.

The angles $\va_{k,i}$, $\va_{k,i}'$ formed by $l_{k,i}$ or
$l_{k,i}'$ with the positive horizontal direction will be called
{\it the inclinations} of $l_{k,i}$, $l_{k,i}'$. Although
$\va_{k,i}=\va'_{k,i}$,  these inclinations will be considered as
independent parameters. It is important to note that each vertex
$A_j'$ of $\hat D$ serves as the origin  for one of the rays
$l_{k,i}$.

Everywhere in this section, $T_1$ will denote the triangle such
that the corresponding ray $l_{1,1}$ has its origin at the vertex
$A_1'$.  Figures~3 and 4 show some notation  used in this section.
\figin{pic5}{700}
 \centerline{Figure 3. Regular proportional
system for small $\theta$}
\vglue12pt
 
For the main parameter $\theta$, we choose the inclination of
$l_{1,1}$: $\theta=\va_{1,1}$. Let $\theta^*$ be the angle formed
by the sides $[A_1',A_2']$ and $[A_1',A_{\hat n}']$ of the convex hull
$\hat D$, then $0<\theta<\theta^*$.

\specialnumber{5}\proclaim{Lemma}  For any  $n$\/{\rm -}\/gon $D$ in standard position there are a finite
number of intervals $(\theta_{j-1},\theta_j)${\rm ,}
$0=\theta_0<\theta_1< \ldots<\theta_{s+1}=\theta^*${\rm ,} such that for
each interval $(\theta_{j-1},\theta_j)$ there is a number $m_j${\rm ,}
 $\hat
n\le m_j\le n$ and a one parameter family of proportional
admissible systems $\{T_i(\theta)\}_{i=1}^{m_j}${\rm ,} which
continuously depend  on~$\theta${\rm ,} $\theta_{j-1}<\theta<\theta_j$
and satisfy the following conditions\/{\rm :}
\begin{itemize}
\ritem{a)} The inclinations  $\va_i(\theta)$ of $l_{1,i}(\theta)${\rm ,}
 $i=1,\ldots,m_j${\rm ,}
strictly increase in $\theta_{j-1}<\theta<\theta_j$.

\ritem{b)} If $\theta\to \theta_j-0$  or $\theta \to \theta_{j-1}+0${\rm ,}
$j=1,\ldots,s${\rm ,} then each triangle $T_i(\theta)$ converges to a
limit triangle $T_i^-(\theta_j)$ or $T_i^+(\theta_{j-1})${\rm ,} some of
which but not all can degenerate to certain nondegenerate
segments. For every $j=1,\ldots,s${\rm ,} the sets of nondegenerate
limit configurations $\{T_i^-(\theta_j)\}$ and
$\{T_i^+(\theta_{j})\}$ coincide.
\end{itemize}

The function $\theta\mapsto\{T_i(\theta)\}_{i=1}^{m(\theta)}$
establishes a one\/{\rm -}\/to\/{\rm -}\/one continuous correspondence between the
interval $(0,\theta^*)$ and the set of all proportional systems
admissible for $D$. \endproclaim

{\it Proof}. 1) First we show that a regular proportional system,
if  it exists for some $\theta^0$, $0<\theta^0<\theta^*$, exists
also for all $\theta$ in some  neighborhood of~$\theta^0$. Let
$\bar\va^0_{\phantom{|}}=(\va_1^0,\ldots,\va_m^0)$ denote the vector of
inclinations for the regular proportional system
$\{T_i^0\}_{i=1}^m$ corresponding to $\theta^0$. Turning all the
rays $l_{1,i}=l_{1,i}(\bar\va^0)$ onto  small angles
$\e_i=\va_i-\va^0_i$ around the origin of $l_{1,i}(\bar\va^0)$, we
get a new system of triangles $\{T_i(\bar\va)\}_{i=1}^m$
corresponding to the vector of inclinations
$\bar\va=(\va_1,\ldots,\va_m)$. This new system is regular and
admissible for $D$, but not proportional in general.

Simple computations show that the coefficient $k_i(\bar \va)$ of
$T_i(\bar \va)$, $i=1,\ldots,m$, is a differentiable function
depending only on $\va_i$ and $\va_{i+1}$ such that\break $\partial
k_i/\partial \va_i<0$, $\partial k_i/\partial \va_{i+1}>0$. Thus
$k_i(\bar\va)$, $i=1,\ldots,m$, satisfy assumptions of Lemma~3.
This implies that for each $\theta=\va_1$ in some small interval
$\theta^0-\delta<\theta<\theta^0+\delta$ there are unique
inclinations $\va_i(\theta)$, $i=2,\ldots,m$ such that the system
$\{T_i(\theta)\}_{i=1}^m$ with $T_i(\theta)=T_i(\bar\va(\theta))$
is regular and proportional.\break\vglue-11pt\noindent Moreover, since
$d\va_i(\theta^0)/d\theta>0$ by inequality (3.4) of Lemma~3, all
rays $l_{1,i}(\theta)$ turn in the same direction as
$l_{1,1}(\theta)$ does.

Let $\theta'<\theta<\theta''$ be the maximal interval containing
$\theta^0$ that carries a regular proportional system. If
$\theta'>0$ then at least one of the $2m$ limit rays
$l_{1,i}(\theta')$, $l_{2,i}(\theta')$, $i=1,\ldots,m$ has two
vertices of $D$. Similarly, if $\theta''<\theta^*$ then at least
one of the $2m$ limit rays $l_{1,i}(\theta'')$,
$l_{2,i}(\theta'')$ has two vertices of $D$. Note that the limit
system of rays always exists since the inclinations
$\va_i(\theta)$ are monotone in $\theta$. In other words, the
limit systems $\{T_i(\theta')\}_{i=1}^m$ and
$\{T_i(\theta'')\}_{i=1}^m$ are singular (=nonregular).

Indeed, if for instance, $\theta'>0$, $\{T_i(\theta')\}_{i=1}^m$
is regular, and all limit triangles $T_i(\theta')$ do not
degenerate, we can continue construction of the system as above
into some right neighborhood of the point $\theta'$.

The degeneracy of $T_i(\theta')$ belonging to the regular limit
system may occur  in two cases. First, if $l_{1,i}(\theta')$ is
parallel to $l_{2,i}(\theta')$: Since the system
$\{T_i(\theta')\}_{i=1}^m$ is regular, $l_{1,i}(\theta')$ and
$l_{2,i}(\theta')$ do not lie on the same straight line. This
implies that $\al_i=\al_i(\theta)\to 0$ and
$\sigma_i=\sigma_i(\theta)\to \infty$ as $\theta\to \theta'+0$.
Since the system $\{T_i(\theta)\}_{i=1}^m$ is proportional for
$\theta'<\theta<\theta''$, the latter implies that
$\al_j(\theta)\to 0$ for all $j=1,\ldots,m$ as $\theta\to \theta'+0$.
This contradicts the condition $\sum_{i=1}^m\al_i(\theta)=2\pi$.

The second type of degeneracy can happen when some triangle
$T_i(\theta')$ shrinks to a point $C\in\partial D$ different from
the vertices of $D$. Since the considered limit system of
triangles is regular, the limit angle $\al_i(\theta')$ cannot be
zero. Hence the limit coefficients
$k_j(\theta')=k_i(\theta')=\infty$. This yields that all the limit
triangles shrink to some points on $\partial D$ different from the
vertices of $D$. This certainly cannot happen, since for every
$i=1,\ldots,\hat n$ and every $\theta$, $A'_i$ belongs to the
boundary of some triangle under consideration.

Thus our analysis show that the limit systems $\{T_i(\theta')\}$
and $\{T_i(\theta'')\}$ are singular. Note that each of them
contains at least $\hat n$ nondegenerate triangles.

\vglue4pt

2) Now we show that a regular proportional system  exists for some
$\theta>0$ small enough. We shall consider rays $l_i=l_i(\e_i)$,
$i=1,\ldots,\hat n$, outgoing from the vertices $A'_i$ of the
convex hull $\hat D$ with inclinations $\va_i=\tilde \va_i+\e_i$,
where $\tilde \va_i$ is the inclination of the side
$[A_i',A'_{i+1}]$ and $0<\e_i\le \e_i^0$. Here  $\e_i^0>0$ are
fixed and small enough such that for $0<\e_i\le \e_i^0$ the ray
$l_i=l_i(\e_i)$ does not contain vertices of $D$ except $A_i'$.

Since $D$ is in general position, it follows that $l_i$ cuts off
from $D$ a triangle with vertices $A'_{i+1}$, $B_i$ and $C'_i$,
where $B_i$ precedes $C'_i$ along $l_i$; see Figure~3.

Let $T_i=T_i(\e_i,\e_{i+1})$ denote the triangle with vertices
$A_{i+1}'$, $B_i$, and $C_i$, where $C_i$ is the point of
intersection of $l_i$ and $l_{i+1}$. Then $\{T_i(\bar
\e)\}_{i=1}^{\hat n}$ with $\bar \e=(\e_1,\ldots,\e_{\hat n})$, is
a regular system admissible for  $D$.
 Let  $\al_i$, and $\sigma_i$ denote
the  base angle and area of $T_i$.  As we observed above, the
coefficient $k_i(\bar\e)=\al_i(\bar \e)/\sigma_i(\bar\e)$ depends
only on $\e_i$ and $\e_{i+1}$. In addition, $k_i(\bar\e)$ is
continuous and strictly increases to $\infty$ as $\e_i$ decreases
to $0$. Therefore the maximal coefficient $$ k=k(\bar\e)=\max_i
k_i(\bar \e)$$ is continuous in $\bar\e$ and $k(\bar\e)\to \infty$
if at least one of the parameters $\e_i$ goes to $0$. This implies
that the minimum
\be
p=\min k(\bar \e) \quad{\mbox{over the set}} \quad 0<\e_i\le
\e_i^0, \quad i=1,\ldots,\hat n,\hskip.5in \ee is achieved at some point
$\bar \e^*=(\e^*_1,\ldots,\e^*_{\hat n})$ with $0<\e_i^*\le
\e_i^0$ for all $i=1,\ldots,\hat n$.

Let us show that $k_i^*=k_i(\bar\e^*)$ equals $p$ for all
$i=1,\ldots, \hat n$. If not,  we consider the set $Q_1(\bar\e^*)$
of vertices $A'_i$ such that $k_i^*=p$ and the set
$Q_2(\bar\e^*)\not=\emptyset$ of vertices $A_i'$ such that
$k_i^*<p$. There is an index $j$ such that $k_j<p$ and
$k_{j+1}=p$. Decreasing $\e_j$ slightly, we get a configuration
$\{T_i(\bar\e)\}_{i=1}^{\hat n}$ for which $k(\bar\e) \le p$ and
the set $Q_1(\bar\e)$ corresponding to this new configuration
contains one vertex less than the \pagebreak set $Q_1(\bar\e^*)$.

If  $Q_1(\bar \e)$ is empty, we get a contradiction to (3.15). If
not, we can repeat the previous procedure. After a finite number
of steps we get a system with $k_i<p$ for all $i=1,\ldots,\hat n$,
contradicting (3.15). This proves that for some $\theta>0$
small enough and, by the same arguments, for some $\theta$ close
enough to~$\theta^*$, there is a regular proportional system.

\vglue4pt

3) For a nonconvex polygon, analysis in 1) and 2) shows that
there are two intervals $(0,\theta_1]$ and $[\theta',\theta^*)$,
each of which carries a parametric family of proportional systems
$\{T_i(\theta)\}_{i=1}^{\hat n}$ with $0<\theta\le \theta_1$ and
$\theta'\le \theta<\theta^*$. Note that for every $\theta$ in
these intervals the system consists of $\hat n$ triangles and the
limit systems $\{T_i(\theta_1)\}_{i=1}^{\hat n}$ and
$\{T_i(\theta')\}_{i=1}^{\hat n}$ are singular.

Next we show that any proportional singular system
$\{T_i(\theta_k)\}_{i=1}^m$ with $0<\theta_k<\theta^*$ and $m$
depending on $k$, can be continued into  some right neighborhood
of $\theta_k$. First we complete the system
$\{T_i(\theta_k)\}_{i=1}^m$ with certain {\it degenerate
triangles} $\hat T_s$ corresponding to singular triangles
$T_{i_s}$ that have two vertices $B_{1,s}$ and $B_{1,s}'$ on the
ray $\tilde l_{1,s}=l_{1,i_s}(\theta_k)$. All possible singular
configurations having this property are depicted in Figure~4 below.
The shaded areas in these figures belong  to the polygon $D$. The
configurations shown in Figure~4: 1, 2, 5, 6, 7, 13, 14, 15, and
16 can also have a second vertex $B_{2,s}'$ on the ray
$l_{2,i_s}$. The angle of $D$ corresponding to this possible
second vertex is shown in the dashed line. Thus, the number of all
possible configurations depicted in Figure~4 equals 26.

For configurations in Figure~4: 15, 16, and 17, $\hat T_s$ will
denote a degenerate triangle having its degenerate base at the
point $C_s\in\partial D$, where the ray $\tilde l_{1,s}$ enters
into $D$ for the first time, and the base vertex at the point
$\hat a_s$ that follows $C_s$ on $\tilde l_{1,s}$ and satisfies
the following condition:
\be
k(\hat T_s):=2|\hat a_s-C_s|^{-2}=k(\theta_k), \ee where
$k(\theta_k)=k_i(\theta_k)$ is the coefficient of the limit system
$\{T_i(\theta_k)\}_{i=1}^m$. 

 For configurations in Figure~4: 1, 2, 3, 4, 12, 13  and 14, we put
$C_s=B_{1,s}'$,  and then   define the point $\hat a_s$, triangle
$\hat T_s$, and the coefficient $k(\hat T_s)$ as above.

Finally, for configurations in Figure~4: 5--11, $\hat T_s$ has
the base $[B_{1,s},B_{1,s}']$ and the base vertex at the point
$\hat a_s$ that follows $B_{1,s}'$ on $\tilde l_{1,s}$ and
satisfies the condition:
\be
k(\hat T_s):=2|(\hat a_s-B_{1,s})(\hat
a_s-B_{1,s}')|^{-1}=k(\theta_k). \ee It is important to note, that
in all cases the point $\hat a_s$, and therefore the triangle  $\hat T_s$, exists and is uniquely determined by  condition (3.16)
or (3.17).

We combine all the triangles $T_i(\theta_k)$ and $\hat T_s$ into a
new system $\{R_i\}_{i=1}^{\tilde m}$,\break $m<\tilde m\le n$, keeping
our usual convention concerning numbering. The latter means, in
particular, that a singular triangle $T_{i_s}$ and the
corresponding degenerate triangle $\hat T_s$ get the indices
$i_s+s-1$ and $i_s+s$ in this new system.

\centerline{\BoxedEPSF{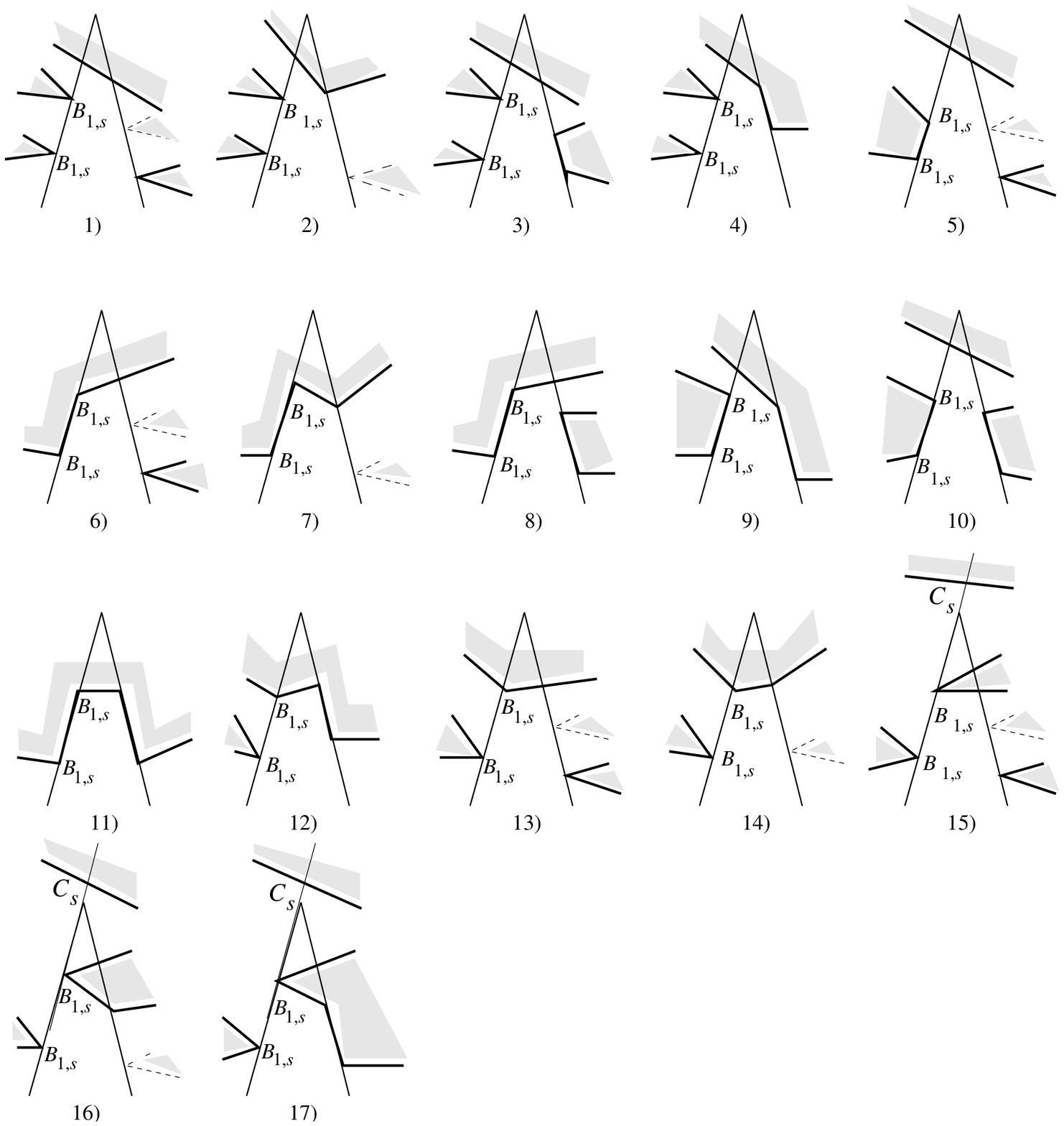 scaled 670}}
 \vglue12pt
\centerline{Figure 4. Singular configurations}
\vglue6pt

To each $R_i$ we associate two rays $r_i$ and $r_{i+1}$ according
to the following rules.  Each regular triangle $T_i$ comes into
the new system with the corresponding rays $l_{1,i}$ and
$l_{2,i}$. If $R_i=T_j$ is singular, then $r_i=l'_{1,j}$,
$r_{i+1}=l_{2,j}$. If $R_i=\hat T_s$ is degenerate corresponding
to a singular triangle $\tilde T_s=T_{i_s}$, then $r_i=l_{1,i_s}$,
$r_{i+1}=l'_{1,i_s}$.

According to our notation, $R_i$ and $R_{i+1}$ have $r_{i+1}$ as a
common associate ray and one can easily check that this is
actually the case.

Let $\bar\e=(\e_1,\ldots,\e_{\tilde m})$, where $\e_i$ is small
enough and  not necessarily positive. We consider a varying system of
rays $r_i(\bar\e)$, $i=1,\ldots,\tilde m$, obtained as follows.

If $r_i$ corresponds to some ray $l_{1,j}$, then $r_i(\bar\e)$ is
obtained from $r_i$ by rotation onto the angle $\e_i$ around the
origin of $r_i$. Thus $r_i(\bar\e)=r_i(\e_i)$ depends only on
$\e_i$ for such rays.

If $r_i$ corresponds to $l'_{1,j}$, then $r_{i-1}$ corresponds to
$l_{1,j}$ and $r_i(\bar\e)$ will denote the ray having a common
origin with $l'_{1,j}$ that intersects $r_{i-1}(\bar\e)$ at the
point $a_{0,i}(\bar\e)$ such that $|a_{o,i}(\bar\e)-B_{1,j}|=|\hat
a_j-B_{1,j}|-\e_i$, where $\hat a_j$ is defined by (3.16) or
(3.17) and $B_{1,j}$ is the origin of $l_{1,j}$. One can easily
see that $r_i(\bar\e)$ corresponding to $l'_{1,j}$ depends on two
parameters: $\e_{i-1}$ and $\e_i$.
\vglue2pt
Let $q_i$ denote the straight line passing through the side of $D$
that contains the base of $R_i$. For configurations shown in
Figure~4: 1, 2, 3, 4, 12, 13, and 14 there are two such straight
lines. In this case, $q_i$ denotes the one through which the
corresponding ray $r_i(\e_i)$ with $\e_i>0$ small enough enters
into $D$ for the first time.

Let $R_i(\bar\e)$ be the triangle (degenerate or not) having its
base on $q_i$ and sides belonging to the rays $r_i(\bar\e)$ and
$r_{i+1}(\bar\e)$. If $R_i(\bar\e)$ is degenerate, we assume in
addition that $R_i(\bar\e)$ has the base vertex at the point
$a_{0,i}(\bar\e)$ such that\break $|a_{0,i}(\bar\e)-B_{1,j}|=|\hat
a_j-B_{1,j}|-\e_i$, where $\hat a_j$ and $B_{1,j}$ are defined
above.

If $\e_i>0$ for all rays $r_i(\e_i)$ depending only on one
parameter, then it is not difficult to see that
$\{R_i(\bar\e)\}_{i=1}^{\tilde m}$ is a regular system admissible
for $D$. Of course, if at least one such $\e_i$ is negative, the
varied system is not admissible. Consider the coefficients
$k_i(\bar\e)$ of the triangles $R_i(\bar\e)$. For the degenerate
triangles, $k_i(\bar\e)$ is defined by (3.16) or (3.17) with $\hat
a_j$ replaced by $a_{0,i}(\bar\e)$ .

By our construction, $k_i(\bar\e)$ depends on two or three
parameters: $\e_i$ and $\e_{i+1}$ or $\e_{i-1}$, $\e_i$,
$\e_{i+1}$, respectively. By direct computation one can easily
check that each $k_i(\bar\e)$ has continuous partial derivatives
near the point $\bar\e^0=(0,\ldots,0)$. Moreover, $\partial
k_i(\bar\e_0)/\partial\e_i<0$, $\partial k_i(\bar\e^0)/\partial
\e_{i+1}>0$ if $k_i$ corresponds to a regular\break\vglue-12pt\noindent triangle; $\partial
k_i(\bar\e^0)/\partial \e_i=0$, $\partial k_i(\bar\e^0)/\partial
\e_{i+1}>0$, if $k_i$ corresponds to a degenerate\break\vglue-12pt\noindent  triangle, and
$\partial k_i(\bar\e^0)/\partial \e_{i-1}<0$, $\partial k_i(\bar
\e^0)/\partial \e_i=0$, $\partial k_i(\bar\e^0)/\partial
\e_{i+1}>0$ if $k_i$ corresponds to a singular triangle.

Therefore, the functions $k_1,\ldots,k_{\tilde m}$ satisfy the
assumptions of Lemma~3. This lemma implies that there is
$\delta>0$ such that for each $\e_1\in (-\delta,\delta)$ there are
unique continuously differentiable functions $\e_i(\e_1)$ solving
the equations \be k_1(\e_1,\ldots,\e_{\tilde
m})=k_2(\e_1,\ldots,\e_{\tilde m} )=\ldots=k_{\tilde m}
(\e_1,\ldots,\e_{\tilde m}) \ee  such that $\e_i(\e_1)\to 0$ as
$\e_1\to 0$.

In addition, inequality (3.4) of Lemma~3 shows that each parameter
$\e_i(\e_1)$ corresponding to some ray $l_{1,j}$ strictly
increases when $\e_1$ does. As we noted above, the latter implies
that the system of triangles $\{R_i(\bar\e(\e_1))\}_{i=1}^{\tilde
m}$ with $\bar\e(\e_1)=(\e_1,\e_2(\e_1),\ldots,\e_{\tilde
m}(\e_1))$ is admissible for $D$. By (3.18), this system is
proportional.

Thus we have proved that every singular proportional system
$\{T_i(\theta_k)\}_{i=1}^m$ can be continued into some right
neighborhood and, by the same arguments, into some left
neighborhood of the parameter $\theta_k$. The arguments above show
also that such continuation is unique. Moreover, since continued
systems are regular, the arguments in the part 1) show that all
inclinations $\va_{k,i}(\theta)$ are monotonic.

4) The arguments in 1)--3) show that there is a family
${\cal{T}}(\theta)$ of proportional, admissible for $D$ systems
$\{T_i(\theta)\}_{i=1}^{m(\theta)}$ that continuously in the sense
of this lemma depend  on the parameter $0<\theta<\theta^*$ and
satisfy conditions a) and b).

To prove the last assertion of the lemma, we assume that there is
a proportionally admissible for $D$ system $\{\tilde
T_i\}_{i=1}^{\tilde m}$ that is not included in
${\cal{T}}(\theta)$. In\break \vglue-11.5pt\noindent this case the arguments above show that
there is a second family $\tilde{\cal{T}}(\theta)$ of systems
$\{\tilde T_i(\theta)\}_{i=1}^{\tilde m(\theta)}$ satisfying the
same assertions of the lemma. The uniqueness of continuation in a
neighborhood established in 1) and 3) shows that
$\{T_i(\theta_1)\}_{i=1}^{m(\theta_1)}\not=\{\tilde
T_i(\theta_2)\}_{i=1}^{\tilde m(\theta_2)}$ for all
$\theta_1,\theta_2\in (0,\theta^*)$.
\vglue2pt
Let $\psi_i(\theta)$ and $\tilde \psi_i(\theta)$ denote the
inclinations of the rays $l_{1,k_i}(\theta)$ and $\tilde
l_{1,s_i}(\theta)$ corresponding to the triangles
$T_{k_i}(\theta)$ and $\tilde T_{s_i}(\theta)$ and outgoing from
the vertex $A'_i$ of $\hat D$. First we show that for every
$i=1,\ldots,\hat n$,
\be
\psi_i(\theta)\to \tilde\va_i\quad {\mbox{and}} \quad
\tilde\psi_i(\theta)\to \tilde\va_i \quad {\mbox{as}}\quad
\theta\to 0,\ee where $\tilde\va_i$ is the inclination of
$[A'_i,A'_{i+1}]$. Suppose for instance, the first relation in
(3.19) is not valid. Then there is an index $j$ such that
\be\psi_j(\theta)\to \tilde\va_j\quad {\mbox{ and}} \quad
\psi_{j-1}(\theta)\to \tilde \va_{j-1}+\e_0 \quad {\mbox{as}}
\quad \theta\to 0,\ee with some $\e_0>0$.

Consider the triangle $T$, possibly infinite, with the base
$[A'_{j-1},A'_j]$ and sides on the limit rays $l_{1,k_j}(0)$,
$l_{1,k_{j-1}}(0)$. (3.20) shows that $T\cap D\not=\emptyset$. Let
$T_k(\theta)$ be a nondegenerate triangle of the system
$\{T_i(0)\}_{i=1}^{m(0)}$ having the associated sector
$S_k(\theta)$ such that $[a_k(\theta),b_k(\theta)]:=\bar
S_k(\theta)\cap [A'_{j-1},A'_j]$ is not empty. Let $R_k(\theta)$
be a triangle with the base $[a_k(\theta),b_k(\theta)]$ that has a
common base angle with $T_k(\theta)$. For $\theta>0$ small enough,
let $\{T_k(\theta)\}_j$ and $\{R_k(\theta)\}_j$ be the systems of
all such triangles $T_k(\theta)$ and $R_k(\theta)$ corresponding
to $[A_{j-1}',A_j]$. Since the inclination $\va_k(\theta)$
corresponding to $T_k(\theta)$ is a monotonic function of $\theta$,
there are limit systems of triangles $\{T_k(0)\}_j$ and
$\{R_k(0)\}_j$ as $\theta\to 0$. Since $D$ is in general position and ${\cal{T}}(\theta)$ consists of proportional systems, (3.20)
implies that $k_i(\theta)\to \infty$ for all $i$ as $\theta\to 0$.
This implies that all limit triangles $T_k(0)$ are degenerate.

On the other hand, since $D$ is in general position and 
$[A'_{j-1},A'_j]$ is covered by the system $\{\bar R_k(0)\}_j$,
the latter system contains nondegenerate triangles. It is clear
that $\{R_k(0)\}_j$ is an admissible system (in the sense of
Lemma~4) for the triangle $T$ defined above. By Lemma~4, $\mathbold{\cup}\bar
R_k(0)\supset T$. Since $T\cap D\not = \emptyset$, the latter
implies that for some $k$, $R_k(0)\cap D\not=\emptyset$. Since
$T_k(0)$ is degenerate, the corresponding limit sector $S_k(0)$
should have a nonempty intersection with $D$. This contradiction
shows that for every $i$, $\psi_i(\theta)\to \tilde \va_i$ and
similarly $\tilde\psi_i(\theta)\to \tilde\va_i$ as $\theta\to 0$.

(3.19) shows that for some $\theta_0>0$ small enough the systems
$\{T_i(\theta_0)\}_{i=1}^{\hat n}$ and $\{\tilde
T_i(\theta_0)\}_{i=1}^{\hat n}$ are of the type considered in part
2) of this proof; cf.\ Figure~3. In particular, each of them
contains $\hat n$ triangles.

Let us show that $\{T_i(\theta_0)\}_{i=1}^{\hat n}=\{\tilde
T_i(\theta_0)\}_{i=1}^{\hat n}$. If not, then there is an index
$i_0$, $1\le i_0\le \hat n-1$, such that $\psi_i(\theta_0)=\tilde
\psi_i(\theta_0)$ for $i=1,\ldots,i_0$ and
$\psi_{i_0+1}(\theta_0)\not=\tilde \psi_{i_0+1}(\theta_0)$. To be
definite, let \be\psi_{i_0+1}(\theta_0)>\tilde
\psi_{i_0+1}(\theta_0).\ee Then \be
k_{i_0}(\psi_{i_0}(\theta_0),\psi_{i_0+1}(\theta_0))>\tilde
k_{i_0}(\psi_{i_0}(\theta_0),\tilde\psi_{i_0+1}(\theta_0)),\ee
where $k_i$ and $\tilde k_i$ denote the coefficients of $T_i$ and
 $\tilde T_i$, respectively. Since the functions $k_i$ and $\tilde k_i$
strictly decrease in their first parameter and strictly increase
in the second one, it follows that $\psi_i(\theta_0)>\tilde
\psi_i(\theta_0)$ for all $i=i_0+1,\ldots,\hat n$. This implies
that $$ k_{\hat n}(\psi_{\hat n}(\theta_0),\theta_0)<\tilde
k_{\hat n}(\tilde\psi_{\hat n}(\theta_0),\theta_0),$$
contradicting (3.22). This contradiction shows that
$\{T_i(\theta_0)\}_{i=1}^{\hat n}=\{\tilde
T_i(\theta_0)\}_{i=1}^{\hat n}$ and therefore
${\cal{T}}(\theta)=\tilde{\cal{T}}(\theta)$ for all $0<\theta<
\theta^*$. This finishes the proof of Lemma~5. \hfill\qed\vglue12pt

Let $D$ be a convex $n$-gon in standard position having the
internal  angle
$\va_i$ at the vertex $A_i$, $i=1,\ldots,n$. Let
$\theta_1=\theta$. For $i=2,\ldots,n$, let
$\theta_i=\theta_i(\theta)$ be functions continuous on $0\le
\theta\le \va_1$ such that $\theta_i(0)=0$,
$\theta_i(\va_1)=\va_i$ and
\be
0<\theta_i(\theta)<\va_i, \quad
\va_{i+1}-\theta_{i+1}(\theta)+\theta_i(\theta)<\pi \quad
{\mbox{for all}}\quad 0<\theta<\va_1.\quad
\ee
Let $l_i(\theta)$ denote
the ray having the angle $\theta_i(\theta)$ with the side
$[A_i,A_{i+1}]$ at $A_i$. If $T_i=T_i(\theta)$ denotes the
triangle with the base $[A_i,A_{i+1}]$ having its sides on the
rays $l_i(\theta)$ and $l_{i+1}(\theta)$, then (3.23) guarantees
that for every $0<\theta<\va_1$ the system
$\{T_i(\theta)\}_{i=1}^n$ is admissible for $D$.

\specialnumber{6}\proclaim{Lemma}  Let $D$ and $\{T_i(\theta)\}_{i=1}^n$ be a convex $n$\/{\rm -}\/gon and
a system of triangles described above. Then there is $\theta^*${\rm ,}
$0<\theta^*<\va_1${\rm ,} such that $\mathbold{\cup}_{i=1}^n \overline
{T}_i(\theta^*)\supset D$. \endproclaim

\demo{Proof} Let $q_i$ denote the straight line containing the
side $[A_i,A_{i+1}]$ and let $d(z,q_i)$ denote the distance from
$z$ to $q_i$. To each $q_i$ we assign the positive weight
$p_i=p_i(\theta)$ as follows. We take $p_0=1$, then for
$i=1,\ldots,n$, we put
\be
p_i=p_{i-1}d(z,q_{i-1})/d(z,q_i)=
p_{i-1}\sin(\va_i-\theta_i(\theta))/
\sin\theta_i(\theta) \quad {\mbox{for}} \quad z\in
l_i(\theta)\setminus A_i.
\ee
Since the quotient
$d(z,q_{i-1})/d(z,q_i)$ is constant on $l_i(\theta) \setminus
A_i$, the functions $p_i(\theta)$ are well defined and continuous
in $0<\theta<\va_1$. Since
$\sin(\va_i-\theta_i(\theta))/\sin\theta_i(\theta) \to +\infty$ as
$\theta\to 0$, $$ p_n(\theta)>\ldots>p_1(\theta)>p_0=1 $$ for all
$\theta$ small enough. Similarly, $$
p_n(\theta)<\ldots<p_1(\theta)<p_0=1 $$ for all $\theta$ close
enough to $\va_1$. Since $p_n(\theta)$ is continuous, there is
$\theta^*$, $0<\theta^*<\va_1$ such that $p_n(\theta^*)=1=p_0$.

We claim that $\{T_i(\theta^*)\}_{i=1}^n$ is a desired cover of
$D$. To show this, we consider components $D_i'$ and $D_i''$ of
$D\setminus l_i(\theta^*)$, where $D_i'$ lies on the left side of
$l_i(\theta^*)$, when we walk on $l_i(\theta^*)$ towards $A_i$.
Using (3.24) one can easily show that for
$i=1,\ldots,n$,
\begin{eqnarray}
d(z,q_{i-1})/d(z,q_i)>p_i(\theta^*)/p_{i-1}(\theta^*)&&{\mbox{for
all}} \quad z\in D_i',\\[5pt]
d(z,q_{i-1})/d(z,q_i)<p_i(\theta^*)/p_{i-1}(\theta^*)&&{\mbox{for
all}} \quad z\in D_i''.
\nonumber
\end{eqnarray}

Assuming that a point $\zeta\in D$ is not covered by $\mathbold{\cup}_{i=1}^n
\overline {T}_i(\theta^*)$, choose an index $j$ such that
$$
p_j(\theta^*)d(\zeta,q_j)=\min_i p_i(\theta^*)d(\zeta,q_i).
$$
Then,
\begin{eqnarray}
d(\zeta,q_j)/d(\zeta,q_{j+1})&\le&
p_{j+1}(\theta^*)/p_j(\theta^*),\\[5pt]
d(\zeta,q_{j-1})/d(\zeta,q_j)&\ge& p_j(\theta^*)/p_{j-1}(\theta^*).
\nonumber
\end{eqnarray}

(3.25) and (3.26) show that $\zeta\in \overline{D}_{j+1}'' \cap
\overline{D}_j'$. Since $\zeta\in D$, the latter implies that
$\zeta\in \overline{T}_j$ contradicting our assumption. \enddemo

{\it Proof of Theorem\/}   4. 1) Assume first that $D$ is in standard
general position. Then by Lemma~5, the set of all proportional
systems admissible for $D$ admits a parametrization in terms of
the angle $\theta$, $0<\theta<\va_1$, formed by the ray
$l_{1,1}(\theta)$ and the segment $[A_1',A_2']$ at the vertex
$A_1'$; see Lemma~5 for the notation. This parametrization of
$\{T_i(\theta)\}$ is continuous in the sense of Lemma~5.

Let $\hat l_k(\theta)$, $k=1,\ldots,\hat n$ be the ray outgoing
from the vertex $A_k'$ of the convex hull $\hat D$ that
corresponds to some triangle $T_j(\theta)$. As mentioned above,
for each $\theta$, every vertex $A_k'$ has one and only one such a
ray.

Let $\hat T_k(\theta)$ denote the triangle with the base
$[A_k',A_{k+1}']$ that has its sides on the rays $\hat
l_k(\theta)$ and $\hat l_{k+1}(\theta)$. It is clear that for
every $\theta$, $0<\theta<\va_1$, $\{\hat
T_k(\theta)\}_{k=1}^{\hat n}$ is a system of triangles admissible
for the convex hull $\hat D$. Moreover, this system is
continuously parametrized by $\theta$, $0<\theta<\va_1$, and
satisfies all other conditions of Lemma~6 for the $\hat n$-gon
$\hat D$. In particular, the limit relations (3.19) are satisfied
as well their counterparts for $\theta\to \va_1$.
Therefore by this lemma, there is $\theta^*$,
$0<\theta^*<\va_1$, such that
\be
\bigcup_{k=1}^{\hat n}\overline{\hat T}_k(\theta^*)\supset \hat D.
\ee Let $S_i(\theta^*)$ and $\hat S_k(\theta^*)$ denote the
sectors associated with the triangles $T_i(\theta^*)$ and $\hat
T_k(\theta^*)$, respectively. Let $I(k)$ denote the set of indices
$i$ such that $S_i(\theta^*)\cap \hat S_k(\theta^*)
\not=\emptyset$. Then $I(1),\ldots,I(\hat n)$ is a disjoint
decomposition of the set of all indices corresponding to the
parameter $\theta^*$.

For $T_i(\theta^*)$ with $i\in I(k)$, let $T_i'(\theta^*)$ be a
triangle with the base on $[A_k',A_{k+1}']$ which has a common
base angle with $T_i(\theta^*)$. The set of the triangles
$\{T_i'(\theta^*):\,i\in I(k)\}$ satisfies the assumptions of
Lemma~4 for the triangle $\hat T_k(\theta^*)$. Therefore,
\be
\bigcup_{i\in I(k)} \overline {T}_i'(\theta^*)\supset \hat
T_k(\theta^*), \quad k=1,\ldots, \hat n. \ee (3.27) and (3.28)
show that
\be
\bigcup_i\overline{T}'_i(\theta^*)\supset \hat D. \ee

Since $T_i(\theta^*)=T_i'(\theta^*)\setminus
\overline{S}_i(\theta^*)$ and $S_i(\theta^*)\cap D=\emptyset$,
(3.29) yields $\mathbold{\cup}_i\overline{T}_i(\theta^*)\break\supset D$. This
proves Theorem~4 for every $n$-gon in general position.

\smallskip

2) For arbitrary $n$-gon $D$ with vertices $A_1,\ldots,A_n$,
consider a sequence of $n$-gons $D^1,D^2,\ldots$, each in general
position, that converges to $D$; i.e., if $D^k$ has vertices
$A_i^k$, $i=1,\ldots,n$, then $A_i^k\to A_i$ as $k\to \infty$. By
part 1), for every $k$ there is a proportional system
$\{T_i^k\}_{i=1}^{m(k)}$ that covers $D^k$. Since $\hat n\le
m(k)\le n$, we may assume that $m(k)=m$ is constant.

Note that
the set of all vertices of all triangles $T_i^k$,$i=1,\ldots,m$,\break 
$k=1,2,\ldots$, is bounded. Therefore we can choose a subsequence
$k_s$, if necessary, such that $T_i^{k_s}$ converge to  limit triangles
$T_i^{\infty}$, $i=1,\ldots,m$,  some of which but not all can degenerate.
It is clear that $\{T_i^\infty\}_{i=1}^m$ is an admissible
proportional system that covers $D$.
\hfill\qed

\section{Proof of Theorem~1}
 
Let $D_n$, $n\ge 3$, be an $n$-gon and let $A=\ar D_n$. By
Theorem~4, there is a proportional admissible for $D_n$ system
$\{T_i\}_{i=1}^m$ with $3\le m\le n$ such that
\be
\sigma:=\sum_{i=1}^m\sigma_i\ge A, \ee where $\sigma_i>0$ denotes
the area of the triangle $T_i$. Let $2\pi\al_i$ be the base angle
of $T_i$. Then
\be
\al_i/\sigma_i=1/\sigma \quad {\mbox{for all}} \quad i=1,\ldots,m,
\ee since $\{T_i\}_{i=1}^m$ is proportional. Let $\{S_i\}_{i=1}^m$
be the system of sectors $S_i$ associated with $T_i$ in the sense
of Section~3. Since $\{T_i\}_{i=1}^m$ is admissible for $D_n$, it
follows that $\{S_i\}_{i=1}^m$ is a competing system of
trilaterals in the sense of Theorem~3 corresponding to a simply
connected domain $\Omega(\overline{D}_n)=\CC
\setminus \overline{D}_n$. Therefore by Theorem~3,
\be
m(\Omega(\overline{D}_n),\infty)\le \sum_{k=1}^m\al_k^2
m(S_k;\infty|a_1^k,a_2^k), \ee where $a_1^k$ and $a_2^k$ are
geometric vertices of $S_k$ different from $\infty$.

By Lemma~1,
\begin{eqnarray}
m(S_k;\infty|a_1^k,a_2^k)&\le&
\frac{1}{2\pi\al_k}\log \frac{4^{\al_k}\al_k
B(1/2,1/2+\al_k)}{(\sigma_k \tan\pi\al_k)^{1/2}}\\[4pt]
&=&
\frac{1}{2\pi\al_k}\log\frac{\pi^{1/2}4^{\al_k}\Gamma(1/2+
\al_k)}{(\sigma_k\tan \pi\al_k)^{1/2}\Gamma(\al_k)}.
\nonumber\end{eqnarray}
Taking into account the proportionality property (4.2) and (4.4),
we get\break from~(4.3)
\be
m(\Omega(\overline{D}_n),\infty)\le \frac{1}{4\pi}\sum_{k=1}^m \al_k\log
\frac{\pi 2^{4\al_k}\Gamma^2(1/2+\al_k)}{\sigma \al_k \tan
\pi\al_k \Gamma^2(\al_k)}= \frac{1}{4\pi}\log\frac{\pi}{\sigma}
+\frac{1}{4\pi}\sum_{k=1}^m H(\al_k),
\ee
where
\be
H(\al)=\al\log\frac{2^{4\al}\Gamma^2(1/2+\al)}{\al \tan\pi\al
\Gamma^2(\al)}.
\ee

In Lemma~7 below we shall show that $H(\al)$ is strictly concave
in $0<\al<1/2$. Since $\sum_{k=1}^m\al_k=1$ and $0<\al_k<1/2$,
(4.5), the concavity property (4.11) and equality (2.14) imply
\begin{eqnarray} \qquad
m(\Omega(\overline{D}_n),\infty)&\le&
\frac{1}{4\pi}\log\frac{\pi}{\sigma}+\frac{m}{4\pi} H(\frac{1}{m})
\\[4pt]
&=&
\frac{1}{4\pi}\log \frac{\pi 2^{4/m} m
\Gamma^2(1/2+1/m)}{\sigma \tan(\pi/m)\Gamma^2(1/m)}
=m(\Omega(\overline{D}^*_m(\sigma)),\infty).
\nonumber
\end{eqnarray}

By Lemma~2, $m(\Omega(\overline{D}_k^*(\sigma)),\infty)$  strictly increases
in $k$ and obviously it strictly decreases in $\sigma$. Therefore,
(4.7) and (4.1) yield
\be
m(\Omega(\overline{D}_n),\infty)\le
m(\Omega(\overline{D}^*_m(\sigma)),\infty)\le
m(\Omega(\overline{D}^*_n(A)),\infty).\ee By (2.2), (4.8) is equivalent to
(1.2).

To prove the uniqueness assertion of Theorem~1, assume that for
$D_n$ considered in the proof above, (1.2) holds with the sign of
equality. By (2.2), the latter is equivalent to the equality for
the reduced modules:
\be
m(\Omega(\overline{D}_n),\infty)=m(\Omega(\overline{D}_n^*(A)),\infty).\pagebreak\ee

In order to have the sign of equality in (4.9),  we must have the
sign of equality in all of the relations (4.1)--(4.8). In
particular, the sign of equality holds in both inequalities in
(4.8), which implies that $m=n$ and $\sigma=A$. The latter
equality shows, in particular, that the triangles $T_i$,
$i=1,\ldots,n$, are mutually disjoint and provide a triangulation
of $D_n$.

Further, equality (4.9) implies that (4.7) holds with the sign of
equality. Since $H(\al)$ is strictly convex, the latter yields
$$
\al_1=\ldots=\al_n=1/n.
$$

To have (4.9), we must have the sign of equality in (4.4) for all
$k=1,\ldots,n$. Now Lemma~1 and (4.2) show that for every
$k=1,\ldots,n$,  $T_k$ is an isosceles triangle having area $A/n$
and the angle $2\pi/n$ at the base vertex $a_0^k$. Therefore, for
every $k=1,\ldots,n$, $S_k$ is an isosceles infinite triangle
having the angle $2\pi/n$ at $\infty$, which is associated with
the triangle $T_k$.

Let $f(\zeta)$ map $\U^*$ conformally onto $\Omega(\overline{D}_n)$ such that
$f(\infty)=\infty$ and let $G_k=f^{-1}(S_k)$. Since all $S_k$ have
the same angle $2\pi/n$ at $\infty$, the uniqueness assertion of
Theorem~3 implies that in the case of equality in (4.3), each
$G_k$ is an infinite sector of the form
$\{\zeta:\,|\zeta|>1,\varphi_k<\arg \zeta
<\varphi_k+2\pi/n\}$. Moreover, geometric vertices of $G_k$
correspond under the mapping $f$ to the geometric vertices of
$S_k$. Now the Schwarz reflection principle implies that $f$ maps
$\U^*$ conformally onto the exterior of a regular $n$-gon. This
finishes the proof of Theorem~1. \hfill\qed\vglue12pt  

 Now we justify the concavity result used in the proof above.

\specialnumber{7}\proclaim{Lemma} 
The function $H(\alpha)$  defined by {\rm (4.6)} is strictly concave
in $0<\alpha<1/2$. In particular{\rm ,}
\be
\sum_{j=1}^n H(\alpha_i)\le
nH(1/n) \quad {\mbox{if}} \quad 0<\alpha_j<1/2 \quad
{\mbox{and}}\quad \sum_{j=1}^n\alpha_j=1.\qquad
\ee
\endproclaim

\demo{Proof} Using the recurrence formula $\al \G(\al)=\G(1+\al)$
and applying the reflection formula $\G(\al)+\G(1-\al)=\pi/\sin
\pi\al$ to the other factor of $F(\al)$ and to one factor of
$\G(1/2+\al)$, we can express $H(\al)$ in a more symmetric form:
$$
H(\al)=4\al^2 \log 2 +\al \log
\frac{\G(1/2+\al)\G(1/2-\al)}{\G(1+\al)\G(1-\al)}.
$$
Differentiating twice, we obtain
\begin{eqnarray*}
H''(\alpha)&=&8\log2+2[\psi(1/2+\alpha)-
\psi(1+\alpha)-\psi(1-\alpha)+\psi(1/2-\al)]\\[4pt]
&&+\ \alpha
[\psi'(1/2+\alpha)-\psi'(1+\alpha)+\psi'(1-\alpha)-
\psi'(1/2-\alpha)].
\end{eqnarray*}

For $\alpha=0$ we use the well-known relations \cite[p.15,18]{BE}
$$
\psi(1)=\g \quad {\mbox{and}} \quad \psi(1/2)=-\g -2\log 2,
$$
where $\g$ denotes the Euler constant \cite[p.1]{BE} to obtain
\be
H''(0)=8\log
2 +4\psi(1/2) -4\psi(1)=0.
\ee

A third differentiation gives: \begin{eqnarray}
\quad\qquad H'''(\alpha)&\hskip-8pt=\hskip-8pt&3\left[\psi'(1/2+\alpha)-\psi'(1+\alpha)+
\psi'(1-\alpha)-\psi'(1/2-\alpha)\right]\\
&\hskip-8pt\hskip-8pt&+\
\alpha\left[\psi''(1/2+\alpha)-\psi''(1+\alpha)-\psi''(1-\alpha)+
\psi''(1/2-\alpha)\right]. \nonumber
\end{eqnarray}

Let $B_1(\alpha)$ and $B_2(\alpha)$ denote expressions in the first
and second brackets in (4.12). Then
\be
B_1(0)=0 \quad
{\mbox{and}}\quad  B_1'(\al)= B_2(\alpha)<0\ee since
$\psi''(t)=-2\sum_{k=0}^\infty (t+k)^{-3}$ increases in $t>0$.
(4.13) shows that $B_1(\al)<0$ and therefore $F'''(\al)<0$ for
$0<\al<1/2$.

The latter and (4.11) imply that $H''(\al)<0$ and therefore
$H(\al)$ is strictly concave in $0<\al<1/2$.

It is well known that concavity of $H$ yields (4.10). \enddemo

\demo{Proof of Lemma~{\rm 2}} The proof follows from Lemma~7 as shown next. Let
$$
\Phi(\al)=\log \frac{2^{4\al} \G^2(1/2+\al)}{\al
\G^2(\al)\tan \pi\al}.
$$
From (2.14) we have
$$
m(\Omega(\overline{D}_n^*(A)),\infty)=(1/4\pi)\,\Phi(1/n)+(1/4\pi)\,
\log (\pi/A).
$$
To show that $\Phi(\al)$ strictly decreases in $\al$, we note that $\Phi$
is given by a difference quotient of the function $H$ in Lemma~7, as
$$
\Phi(\al)=\frac{H(\al)-H(0)}{\al},
$$
since $H(0)=0$. This difference quotient is a strictly decreasing
function of $\al\in (0,1/2)$, by the concavity of $H$. This proves
Lemma~2.
\enddemo

\end{document}